\documentclass[10pt]{amsart}

\usepackage{amsmath, latexsym, amssymb, verbatim, graphics, textcomp}
\usepackage[all,2cell,dvips]{xy}

\theoremstyle{plain}

\newtheorem{thm}{Theorem}
\newtheorem{lem}[thm]{Lemma}
\newtheorem{prop}[thm]{Proposition}
\newtheorem{fact}[thm]{Fact}
\newtheorem{cor}[thm]{Corollary}

\newcommand{\nc}{\newcommand}
\nc{\dmo}{\DeclareMathOperator}

\dmo{\rk}{rk}
\dmo{\Comm}{Comm}
\dmo{\Aut}{Aut}
\dmo{\Out}{Out}
\dmo{\Homeo}{Homeo}
\dmo{\T}{Tv}

\dmo{\Mod}{Mod}
\dmo{\PMod}{PMod}

\nc{\A}{{\mathcal A}}
\nc{\C}{{\mathcal C}}
\renewcommand{\P}{\mathcal{P}}

\nc{\ct}{\widetilde{C}}
\nc{\at}{\widetilde{A}}

\nc{\Cn}{B}
\nc{\cn}{A(\Cn_n)}
\nc{\cnt}{A(\ct_{n-1})}
\nc{\ant}{A(\at_{n-1})}
\nc{\an}{A(A_n)}
\nc{\pa}{P(A_n)}

\nc{\Q}{{\mathbb Q}}
\nc{\Z}{{\mathbb Z}}
\nc{\R}{{\mathbb R}}

\nc{\bpf}{\begin{proof}}
\nc{\epf}{\end{proof}}

\nc{\set}[1]{\{#1\}}

\nc{\p}[1]{{\bf #1} }

\nc{\pics}[3]{\epsfysize=#3 cm \begin{figure}[htb]
\center{\leavevmode \epsfbox{#1.eps}} \caption{#2}  \label{#1pic}
\end{figure}}

\nc{\bl}{ \begin{list}{$\cdot$}{
\setlength{\leftmargin}{.25in}
\setlength{\rightmargin}{.5in}
\setlength{\parsep}{0.5ex plus .2ex minus 0ex}
\setlength{\itemsep}{0.2ex plus 0.2ex minus 0ex}
}
}

\nc{\el}{\end{list}}

\nc{\margin}[1]{\marginpar{\scriptsize #1}}

\begin{document}

\input{epsf.sty}

\title{\bf Injections of Artin groups}

\author{Robert W. Bell}
\author{Dan Margalit}

\address{Department of Mathematics\\ University of Utah\\ 155 S 1440 East \\ Salt Lake City, UT 84112-0090}

\thanks{This first author is partially supported by a VIGRE postdoctoral position under NSF grant number 0091675 to the University of Utah.  The second author is supported by an NSF postdoctoral fellowship.}

\email{rbell@math.utah.edu, margalit@math.utah.edu}



\maketitle

\begin{center}\today\end{center}

\begin{abstract} We study those Artin groups which, modulo their
centers, are finite index subgroups of the mapping class group of a
sphere with at least 5 punctures.  In particular, we show that any
injective homomorphism between these groups is parameterized by a
homeomorphism of a punctured sphere together with a map to the
integers.
We also give a generating set for the automorphism group of the pure braid
group on at least 4 strands.  The technique, following Ivanov, is to
prove that every superinjective map of the complex of curves of a
sphere with at least 5 punctures is induced by a
homeomorphism.\end{abstract}


\section{Introduction}

We investigate injective homomorphisms between Artin groups which,
modulo their center, embed as finite index subgroups in the mapping class group of an $m$-times punctured sphere $S_m$, where $m \geq 5$.

The \emph{extended mapping class group} of a
surface $F$ is the group of isotopy classes of homeomorphisms of $F$:
\[\Mod(F) = \pi_0(\Homeo(F)). \]

\begin{thm} \label{main1} Let $m \geq 5$.
If $G$ is a finite index subgroup of $\Mod(S_m)$ and $\rho:G \to \Mod(S_m)$
is an injective homomorphism, then there is a unique $f \in \Mod(S_m)$ so that
$\rho(g) = fgf^{-1}$ for all $g \in G$.
\end{thm}

In particular, the Theorem~\ref{main1} applies to four infinite families of
Artin groups modulo their centers: $\an/Z$, $\cn/Z$, $\cnt$, and
$\ant$ where $n=m-2$ (see below for definitions).  Throughout, $Z$ denotes the
center of the ambient group; the groups $\cnt$ and $\ant$ have trivial
center.  Thus, Theorem~\ref{main1} is a generalization of work of Charney--Crisp,
who computed the automorphism groups of the aforementioned Artin
groups using similar techniques \cite{cc}.

Following Ivanov, we prove Theorem~\ref{main1} by translating the problem into one about the
\emph{curve complex} $\C(S_m)$.  This is the abstract simplicial flag complex with vertices corresponding to isotopy classes of essential curves in $S_m$ and edges corresponding to disjoint pairs of curves.  To this end, we focus on particular elements of $G$: powers of Dehn twists; each such
element is associated to a unique isotopy class of curves in $S_m$
(see Section~\ref{bg}).  We
show that the injection $\rho$ must take powers of Dehn twists to
powers of Dehn twists, thus giving an action $\rho_\star$ on the vertices of
$\C(S_m)$.  Since $\rho_\star$ is easily seen to be superinjective in the
sense of Irmak (i.e. $\rho_\star$ preserves disjointness and
nondisjointness), we will be able to derive Theorem~\ref{main1} from the
following theorem.

\begin{thm} \label{superinj}
Let $m \geq 5$.  Every superinjective map of $\C(S_m)$ is
  induced by a unique element of $\Mod(S_m)$.
\end{thm}

The proofs of both theorems are modeled on previous work of Ivanov,
who showed that every isomorphism between finite index subgroups of
$\Mod(F)$ is the restriction of an inner automorphism of $\Mod(F)$,
when the genus of $F$ is at least 2 \cite{ni}.  To do this, he
applied his theorem that every automorphism of  $\C(F)$ is induced by an element of $\Mod(F)$.
His method has been used to prove similar theorems by various other authors
\cite{mk,fl,im,ei,ei2,ei3,dm,iim,fi,mv,brm}. In particular, Korkmaz proved that every automorphism of
$\C(S_m)$ is induced by an element of $\Mod(S_m)$ \cite{mk},
and Irmak showed that every superinjective map of $\C(F)$, for
higher genus $F$, is induced by an element of $\Mod(F)$, thus
obtaining the analog of Theorem~\ref{main1} for surfaces of genus at
least 2 \cite{ei,ei2,ei3}.

After the completion of the work presented in this paper, the final cases of Irmak's theorem were completed by Behrstock--Margalit~\cite{bem}, Shackleton~\cite{kjs}, and a subsequent paper by the authors \cite{si}.

\p{Artin groups.} Before we explain the applications of Theorem~\ref{main1} to Artin groups, we recall the basic definitions.  An
\emph{Artin group} is any group with a finite set of generators
$\{s_1, \dots, s_n\}$ and, for each $i \neq j$, a defining relation of
the form
\[s_i s_j \cdots = s_j s_i \cdots \]
where $s_i s_j \cdots$ denotes an alternating string of $m_{ij} =
m_{ji}$ letters. The value of $m_{ij}$ must lie in the set $\{2, 3,
\dots, \infty\}$ with $m_{ij}= \infty$ signifying that there is no
defining relation between $s_i$ and $s_j$.

\begin{figure}[htb]
\begin{center}
\[
\begin{picture}(0,0)(0,0)
\put(185,108){$A_n$}
\put(185,81){$\Cn_n$}
\put(185,53){$\ct_{n-1}$}
\put(185,2){$\at_{n-1}$}
\put(215,108){$(n > 0)$}
\put(215,81){$(n > 1)$}
\put(215,53){$(n > 3)$}
\put(215,2){$(n > 2)$}
\end{picture}
\scalebox{.75}{\includegraphics{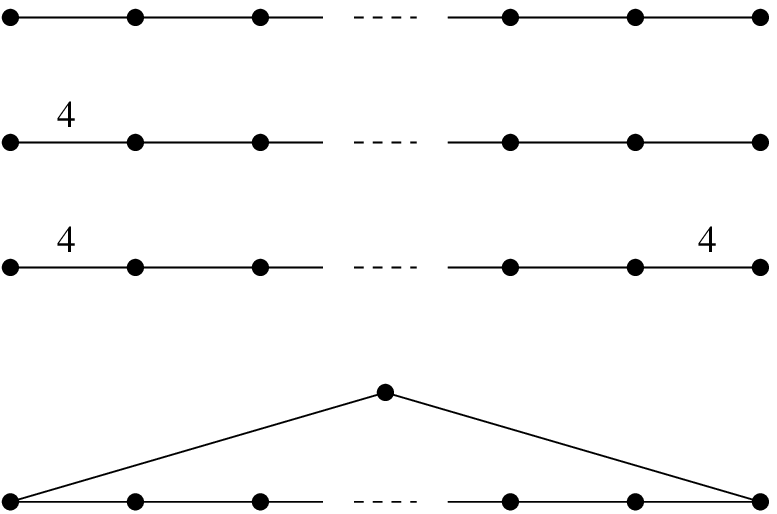}}
\]
\caption{Coxeter graphs with $n$ vertices.}
\label{coxgraphpic}
\end{center}
\end{figure}

It is convenient to define an Artin group by a {\em Coxeter graph},
which has a vertex for each generator $s_i$ and an edge labelled
$m_{ij}$ connecting the vertices corresponding to $s_i$ and $s_j$ if
$m_{ij} > 2$.  The label 3 is suppressed.  The Coxeter graphs $A_n$,
$\Cn_n$, $\ct_{n-1}$, and $\at_{n-1}$ for the Artin groups $\an$,
$\cn$, $\cnt$, and $\ant$ are displayed in Figure~\ref{coxgraphpic}.

\p{Artin groups and mapping class groups.} The Artin group $\an$ is better
known as the braid group on $n+1$ strands. If $D_{n+1}$ is the disk with $n+1$
punctures, and $\Homeo ^+(D_{n+1},\partial D_{n+1})$ is the space of (orientation preserving)
homeomorphisms of $D_{n+1}$ which are the identity on the boundary, then $\an$ is isomorphic to
\[ \Mod(D_{n+1},\partial D_{n+1}) = \pi_0(\Homeo ^+(D_{n+1},\partial D_{n+1})) \]
(see, e.g., \cite{bb}).  The \emph{pure braid group} $\pa$ is the (finite index)
subgroup of $\an$ consisting of elements which fix each puncture of $D_{n+1}$.
The group $\cn$ is isomorphic to a subgroup of $\an$ fixing one given
puncture (see~\cite{da} or~\cite{cc}).

The center of $\an$ is generated by the Dehn twist
about a curve isotopic to $\partial D_{n+1}$; we denote this element by $z$.  Both $\cn$ and $\pa$ inherit the
same center.

We can also identify $\an/Z$, $\cn/Z$, and $\cnt$ with the subgroups
of $\Mod(S_{n+2})$ consisting of orientation preserving elements which
fix one, two, and three particular punctures, respectively; further
$\pa/Z$ is isomorphic to the \emph{pure mapping class group}
$\PMod(S_{n+2})$, which is the finite index subgroup of
$\Mod(S_{n+2})$ consisting of orientation preserving elements which
fix every puncture.  The group $\ant$ is also isomorphic to a finite
index subgroup of $\Mod(S_{n+2})$.  A complete description of these
isomorphisms appears in the paper of Charney--Crisp~\cite{cc}.  The
proofs are due to Allcock, Kent--Peifer, Charney--Peifer, Crisp, and
Charney--Crisp \cite{da,cc,cp,jc,kp}.

\p{Applications.} We now give some consequences of Theorem~\ref{main1}.  A group is \emph{co-Hopfian}
if each of its injective endomorphisms is an isomorphism.

\begin{cor} \label{cohopf}
For $n \geq 3$, all finite index subgroups of $\Mod(S_{n+2})$ are co-Hopfian;
in particular, the groups $\an/Z$, $\cn/Z$, $\cnt$, $\ant$, and
$\pa/Z$ are co-Hopfian.\end{cor}

For each $0 \leq k \leq m = n+2$, let $G_k$ be a subgroup of
$\Mod(S_{m})$ consisting of orientation preserving elements which
fix $k$ given punctures.  Note that $G_0$ is the index 2 subgroup of
$\Mod(S_{m})$ consisting of orientation preserving elements,
$G_1 \cong \an/Z$, and $G_{m-1} = G_{m} = \PMod(S_{m})$.  Also, $G_2
\cong \cn/Z$ and $G_3 \cong \cnt$.

\begin{cor} \label{injgraph} Suppose $n \geq 3$ and let $G$ and $H$ be any of
the groups in Figure~\ref{graphpic}.  Then there exists an injection
$\rho:G \to H$ if and only if there is directed path from $G$ to $H$
in Figure~\ref{graphpic}.
\end{cor}

A concrete way to verify this corollary is to compare the indices of
the groups $G$ and $H$ in $\Mod(S_m)$.

One might also ask whether or not any injection from
Corollary~\ref{injgraph} is unique up to automorphisms of $H$.  The
answer is no.  For instance, since $G_m$ is normal in $G_0$, we may
conjugate $G_m$ by any element of $\Mod(S_m)$ to get an injective
homomorphism $G_m \to G_k$ for any $k$.  However, if $k > 0$, then
$f$ might not fix the $k$ punctures fixed by $G_k$, and so there is
no automorphism of $G_k$ which achieves the injection.

\begin{figure}
\scalebox{.75}{ \xymatrix{ & & & & & \ant \ar@{->}[d] \\
\PMod(S_{n+2}) \ar@{->}[r] & G_{n} \ar@{->}[r] & \dots \ar@{->}[r]&
G_4 \ar@{->}[r] & \cnt \ar@{->}[r] & \cn/Z \ar@{->}[r] & \an/Z
\ar@{->}[r] & G_0 \\} } \caption{The diagram for
Corollary~\ref{injgraph}.} \label{graphpic}
\end{figure}

We are also able to characterize injections between the groups
$\an$, $\cn$, and $\pa$ (with their centers).  There are inclusions:
$\pa \to \cn \to \an$ (see Section~\ref{cat}); all other injections
between these groups are described by the following corollary to the
Theorem~\ref{main1}.

\begin{thm} \label{maininj}
Suppose $n \geq 3$.  Let $G$ be a finite index subgroup of $\an$. If
$\rho: G \to \an$ is an injective homomorphism, then there is an induced
injection $G/Z \to \an/Z$.  Moreover, there is a unique $f \in
\Mod(S_{n+2})$ so that, after identifying $\an/Z$ with the group
$G_1$, we have
\[ \rho(g) Z = f (gZ) f^{-1}\]
for all $g \in G$.
\end{thm}

In Section~\ref{cat}, we explain how this theorem may be applied to
give an explicit list of all injections of $\an$, $\cn$, and $\pa$
into $\an$.  The case of $\an$ was already handled in a previous
paper of the authors \cite{bm}.

Combining Theorem~\ref{maininj} with Corollary~\ref{injgraph}, we immediately obtain an analogue of Corollary~\ref{injgraph} for $\an$.  Precisely, if $L_k$ is a subgroup of $\an$ corresponding to elements which fix $k$ particular punctures, then there is an injective homomorphism $L_j \to L_k$ if and only if $j>k$.  In particular, there is an injective homomorphism between two of the groups $\an$, $\cn$, and $\pa$ if and only if there is an obvious one.

Finally, combining a theorem of Korkmaz with our
understanding of injections of $\pa$ into $\an$, we will prove the
following two results.  The first is a theorem of Charney--Crisp.

\begin{thm} \label{autcn}
For $n \geq 3$, we have
\[ \Aut(\cn) \cong (\Z_2 \times \Z_2) \ltimes (G_2 \times \Z) \]
where $G_2$ is the group defined above.
\end{thm}

Whereas the Charney--Crisp proof of Theorem~\ref{autcn} relies on a semidirect product decomposition of $\cn$ due to Kent--Peifer, we work directly from the isomorphism $\cn \cong L_1$.

\begin{thm} \label{autpa}
Suppose $n \geq 3$ and let $N = {n+1 \choose 2}$.  We have a short exact sequence
\[ 1 \to \Z_2 \ltimes \Z^{N-1} \to \Aut(\pa) \to \Aut(\pa/Z) \to 1 \]
where the third map is the natural one.
\end{thm}

A consequence of Theorem~\ref{autpa} is that we get a generating set for $\Aut(\pa)$ from the standard generating sets for $\Z_2 \ltimes \Z^{N-1}$ and $\Aut(\pa/Z) \cong \Mod(S_{n+2})$ (the last isomorphism is a theorem of Korkmaz).

\p{Acknowledgements.}  The authors would like to thank Benson Farb for
bringing the work of Charney--Crisp to our attention and for his
continued encouragement on this project.  Joan Birman, Chris
Leininger, and Luis Paris were very generous with their time and
energy; for this we are thankful.  We are also grateful to Mladen
Bestvina, Tara Brendle, Ruth Charney, Andy Putman, Michah Sageev, and Steven Spallone
for helpful conversations.  Finally, we thank the referee for comments
which improved the paper.


\section{Background}\label{bg}

\p{Curves.} By a {\em curve} in a surface $F$, we mean the isotopy
class of a simple closed curve in $F$ which is not isotopic to a
point, a puncture, or a boundary component of $F$.  We will often not
make the distinction between a representative curve and its isotopy
class.

We denote by $i(a,b)$ the {\em geometric intersection number}
between two curves $a$ and $b$.

A maximal collection of pairwise disjoint curves in $S_m$ is called a
\emph{pants decomposition}.  Any pants decomposition of $S_m$ or
$D_{m-1}$ has $m-3$ curves.

A curve in $D_{n+1}$ with $k$ punctures in its interior is called a \emph{$k$-curve}.

\p{Curve complex.} The \emph{curve complex} $\C(F)$ for a surface $F$,
defined by Harvey, is the abstract simplicial flag complex with a
vertex for each curve in $F$ and edges corresponding to
geometric intersection zero \cite{wjh}.

A map $\phi:\C(F) \to \C(F)$ is
called \emph{superinjective} if for any two vertices $v$ and $w$ of $\C(F)$, thought of as curves in $F$, we have $i(v,w)=0$ if and only if
$i(\phi(v),\phi(w))=0$.  Superinjective maps of $\C(S_m)$ are injective
for $m \geq 5$ since, given two distinct curves, there is a curve which is disjoint from one but not the other.

\p{Twists.} A \emph{Dehn twist} about a curve $a$, denoted $T_a$ is the element of the mapping class group which has support on an annular neighborhood of $a$,
and is described on that annulus by Figure~\ref{dtpic}.

\pics{dt}{Dehn twist about a curve $a$.}{4}

If $a$ is a 2-curve, we define the \emph{half twist}
about $a$, denoted $H_a$, to be the element of the mapping class group
which has support the interior of $a$, and is described inside this twice-punctured disk by Figure~\ref{htpic}.

\pics{ht}{Half twist about a curve $a$.}{2.5}

For each $f \in \Mod(S_m)$, let $\epsilon(f) = 1$ if $f$ preserves
orientation and $\epsilon(f) = -1$ if not.  We will use the following
connection between the topology and algebra of Dehn twists in $\Mod(S_m)$.

\begin{fact} \label{F1}
Suppose $f \in \Mod(S_m)$.  Then $fT_af^{-1} = T_{f(a)}^{\epsilon(f)}$.
In particular, $[f,T_a] = 1$ implies $f(a) = a$, and powers of Dehn
twists commute if and only if the curves have geometric intersection zero.
\end{fact}

For a group $\Gamma$, we define its {\em rank}, $\rk(\Gamma)$, to be the
maximal rank of a free abelian subgroup of $\Gamma$.  It follows from
work of Birman--Lubotzky--McCarthy that for any surface $F$, $\rk
\Mod(F)$ is realized by any subgroup generated by powers of Dehn twists
about curves forming a pants decomposition for $F$ \cite{blm}; thus,
$\rk \Mod(S_m) = m-3$.  The following theorem of Ivanov gives another
connection between the algebra and topology of $\Mod(S_m)$ \cite{ni}.  We restrict our attention here to the genus 0 case, which has a particularly simple statement.

\begin{thm}\label{T:ivanov}
Let $m \geq 5$ and let $P$ be a finite index subgroup of $\PMod(S_m)$.  An element $g$ of $P$
is power of Dehn twist if and only if $Z(C_P(g)) \cong \Z$ and
$\rk C_P(g) = m-3$.
\end{thm}

We now state a group theoretical lemma, due to Ivanov--McCarthy \cite{im},
which will be used in Proposition~\ref{P:almost}.

\begin{lem}\label{L:gptheory} Let $\rho: \Gamma \to \Gamma'$ be any injective homomorphism
of groups, where $\rk \Gamma' = \rk \Gamma < \infty$.
If $G < \Gamma$ is a free abelian subgroup of maximal rank, and $g \in G$.
Then \[\rk Z(C_{\Gamma'}(\rho(g))) \leq \rk Z(C_{\Gamma}(g)).\]
\end{lem}

Note that Lemma~\ref{L:gptheory} applies whenever $g$ is a power
of a Dehn twist and both $\Gamma$ and $\Gamma'$ are finite index subgroups of
$\PMod(S_m)$.



\section{Subgroups of $\Mod(S_m)$}
\label{pmt}

Let $\rho:G \to \Mod(S_m)$ be an injective homomorphism, where $G$ is a finite index
subgroup of $\Mod(S_m)$ with $m \geq 5$.

\begin{prop} \label{P:almost}
For each curve $a$ in $S_m$, there are nonzero integers $k$ and $k'$ and
a curve $a'$ such that $\rho(T^k_a) = T^{k'}_{a'}$.
\end{prop}

\bpf Let $Q = \PMod(S_m)$, and let $P = Q \cap \rho^{-1}(Q)$.  Since $P$ is a
finite index subgroup of $\Mod(S_m)$, we can choose a $k$ so that $g=T_a^k$
belongs to $P$.  By Theorem~\ref{T:ivanov}, $Z(C_P(g)) \cong \Z$.
Lemma~\ref{L:gptheory} and the fact that $\rho$ is injective imply
that $Z(C_{Q}(\rho(g)) \cong \Z$. Since $\rk \rho(C_P(g)) =
\rk \Mod(S_m)$, Theorem~\ref{T:ivanov} says that $\rho(g)$ must be a
power of a Dehn twist.
\epf

By Proposition~\ref{P:almost}, $\rho$ induces a well-defined action
$\rho_\star$ on curves given by
\[ \rho(T^k_a) = T^{k'}_{\rho_\star(a)}. \]
Applying Fact~\ref{F1}, we have:

\begin{prop} \label{rssi}
The map $\rho_\star$ is a superinjective map of $\C(S_m)$.\end{prop}

We are now ready to complete the proof of Theorem~\ref{main1}, assuming
Theorem~\ref{superinj}.

\bpf By Propositions~\ref{P:almost} and \ref{rssi}, the injection $\rho$
gives rise to a superinjective map $\rho_\star$ of $\C(S)$, which by Theorem~\ref{superinj}
is induced by a unique $f \in \Mod(S_m)$; that is to say, $\rho_\star(c) = f(c)$
for every curve $c$.  Since $f$ is unique, we can check that $\rho(g) = fgf^{-1}$ by checking that $fg(c)=\rho(g)f(c)$ for any curve $c$.
\[T_{fg(c)}^{k'} = \rho(T_{g(c)}^k) = \rho(gT_c^{\pm k}g^{-1}) =\]
\[\rho(g)\rho(T_c^{\pm k})\rho(g)^{-1} = \rho(g)T_{f(c)}^{k''}\rho(g)^{-1} =
T_{\rho(g)f(c)}^{\pm k''}\] Thus, $T_{fg(c)}^{k'} = T_{\rho(g)f(c)}^{\pm k''}$,
which implies that $fg(c) = \rho(g)f(c)$.
\epf


\section{Subgroups of $\an$}

Let $G < \an$ be a finite index subgroup and $\rho:G \to
\an$ an injective homomorphism. To prove Theorem~\ref{maininj}, we
need to show that $\rho$ induces an injective homomorphism $G/Z \to
\an/Z$ and apply Theorem~\ref{main1}.

As with Theorem~\ref{main1}, we shall require the existence of a
superinjective map $\rho_\star$ of $\C(D_{n+1})$ which is induced by $\rho$ in
the sense that for any curve $a$ we have
\[ \rho(T_a^k) = T_{\rho_\star(a)}^{k'} z^{k''} \]
for some integers $k$, $k'$, and $k''$ ($k$ and $k'$ nonzero); as
usual $z$ is the generator of the center of $\an$.  The
argument is exactly the same as in Proposition~\ref{P:almost}, with
Theorem~\ref{T:ivanov} replaced by the following corollary of
Theorem~\ref{T:ivanov}.

\begin{cor} Let $P$ be a finite index subgroup of $\pa$.  An element
  $g$ of $P$ is the product of a central element and a nontrivial power of a
  noncentral Dehn twist if and only if $Z(C_P(g)) \cong \Z^2$
  and $\rk C_P(g)=n$.\end{cor}

We now prove the theorem.

\begin{proof}[Proof of Theorem~\ref{maininj}]

Let $G$ be a finite index subgroup of $\an$ and $\rho:G \to \an$ an
injective homomorphism. We know that $G$ has nontrivial center
$Z(G)$ since it is finite index in $\an$. Further we have $Z(G)
\cong \Z$. Indeed, if $\zeta$ is an element of $Z(G)$, then $\zeta$
must fix every curve in $D_{n+1}$ by Fact~\ref{F1} and the fact that
$G$ is finite index; hence $\zeta$ is a power of $z$.

Let $\zeta$ denote a generator of $Z(G)$.  We now show that
$\rho(Z(G)) < Z$ by showing $\rho(\zeta) \in Z$.  Since $\rk G = \rk
\an$, we have that $\rho(\zeta^k) \in Z$ for some nonzero $k$.

Choose a pants decomposition $\P$ of $D_{n+1}$.  As in
Section~\ref{pmt}, we know that $\rho_\star(\P)$ is also a pants
decomposition ($\rho_\star$ exists by the discussion at the start of this
section). Further, because $\zeta$ is central and $\rho$ is
injective, it follows that $\rho(\zeta)$ fixes each element of
$\rho_\star(\P)$. Since an orientation preserving element of $\Mod(S_3)$ is
determined by its action on the punctures, it follows that
$\rho(\zeta)$ lies in the free abelian subgroup generated by
half twists and Dehn twists in the curves of $\rho_\star(\P)$. Since
$\rho(\zeta)^k \in Z$, it now follows that $\rho(\zeta) \in Z$.

\nc{\rb}{\bar{\rho}}

Moreover, we have that $\rho^{-1}(Z) < Z(G)$, by the injectivity of
$\rho$.  Thus, $\rho$ induces a well-defined injection $G/Z(G) \to
\an/Z$.  Since $G/Z(G)$ is finite index in $\an/Z$, we may apply Theorem~\ref{main1}.  Thus, fixing an identification $\an/Z <
\Mod(S_{n+2})$, there is a unique $f \in \Mod(S_{n+2})$ so that
\[ \rho(g) Z = f(gZ)f^{-1} \]
for all $g \in G$.  This proves the theorem.
\end{proof}

We now take a moment to interpret Theorem~\ref{maininj} in a way that
will be useful to us in the next section.  The element $f \in
\Mod(S_{n+2})$ from the theorem does not necessarily correspond to an
element of $\Mod(D_{n+1})$, for it may switch the puncture in
$S_{n+2}$ corresponding to the boundary of $D_{n+1}$ with another
puncture.  However, even in this case, since $\C(D_{n+1}) \cong
\C(S_{n+2})$, the element $f$ induces an automorphism $f_\star$ of
$\C(D_{n+1})$.  Since $f$ is an element of $\Mod(S_{n+2})$ as opposed
to $\Mod(D_{n+1})$, the map $f_\star$ may take a $k$-curve to an $(n-k+2)$-curve.

Now, let $g$ be a power of a
noncentral Dehn twist or half twist in $G$; for concreteness,
$g=T_a^k$. The coset $gZ$ is a power of a Dehn twist, also denoted
$T_a^k$, thought of as an element of $\Mod(S_{n+2})$. The conjugate
$f(gZ)f^{-1}$ is equal to $T_{f_\star(a)}^{\epsilon(f) k}$. It follows from Theorem~\ref{maininj} that
$\rho(g)$ is a product of $T_{f_\star(a)}^{\epsilon(f) k}$ with a central element.  We again emphasize that $f$ is an element of the mapping
class group of $S_{n+2}$, and not $D_{n+1}$, and so $f_\star$ can
take a curve to one which is not topologically equivalent in $D_{n+1}$.  In fact, we will see
examples of this ``nongeometric'' phenomenon in the next section,
where we classify injective homomorphisms of $\an$, $\cn$, and $\pa$
into $\an$.  As each of these groups is generated by half twists and
Dehn twists, we will be able to understand these injections via the following corollary to Theorem~\ref{maininj}, which summarizes the above discussion.

\begin{cor}
\label{twist image}
Let $G$ be a finite index subgroup of $\an$ and $\rho:G \to
\an$ an injective endomorphism.  There is an $f \in \Mod(S_{n+2})$
so that for any power of a Dehn twist $T_a^k \in G$, we have
\[ \rho(T_a^k) = T_{f_\star(a)}^{\epsilon(f)k} z^{t} \]
for some integer $t=t(T_a^k)$.
\end{cor}

The analogous statement for half twists also holds.  We remark that the reason we focus on Dehn twists and half twists here is that in $\an$ there is a natural representative of a Dehn twist coset of $\an/Z$, and so, combined with the action of $f_\star$ on curves, there is a relatively simple form for the image under $\rho$ of a power of a twist.

\p{Moving punctures criterion.} Another fact which will be useful in
the next section is that $f$ must send moving punctures to moving
punctures; that is, the set of punctures of $S_{n+2}$ which are not
fixed by every element of $G/Z(G)$ must be sent by $f$ to into
the $n+1$ punctures which are not fixed by $\an/Z$ (recall that only
one puncture is fixed by $\an/Z$).  This is because
conjugation by $f$ sends fixed punctures to fixed punctures and
moving punctures to moving punctures.  Below, we call this the
\emph{moving punctures criterion}.  We remark that this criterion can be used to derive Corollary~\ref{injgraph} from Theorem~\ref{main1}.


\section{Catalogue of injections}\label{cat}

We now use Theorem~\ref{maininj} to list all injections of the groups
$\pa$, $\cn$, and $\an$ into $\an$.  As usual, we denote the generator
of $Z$ by $z$.

Instead of applying Theorem 4 directly, we will instead use Corollary~\ref{twist image} and the moving punctures criterion.  We use the notation of Corollary~\ref{twist image}: given an element $f$ of $\Mod(S_{n+2})$, the symbol $f_\star$ denotes the induced automorphism of $\C(S_{n+2})$; the identification $\C(D_{n+1}) \cong \C(S_{n+2})$ comes from the identification $\an/Z < \Mod(S_{n+2})$.

\p{Injections of \boldmath$\an$.} The Artin group $\an$ is defined via the presentation given by Figure~\ref{coxgraphpic}.  We denote the generators by $\sigma_1, \dots, \sigma_n$.  Under the identification with $\Mod(D_{n+1},\partial
D_{n+1})$, each generator $\sigma_i$ corresponds to a half twist $H_{a_i}$ about a curve $a_i$ in $D_{n+1}$ (see \cite{bb}).

Let $\rho:\an \to \an$ be an injective homomorphism.  Applying the
moving punctures criterion, we see that the element $f \in \Mod(S_{n+2})$ given
by Theorem~\ref{maininj} must send the puncture of $S_{n+2}$ fixed
by $\an/Z$ to itself.  Therefore, we may think of $f$ as an element
of $\Mod(D_{n+1})$.  Corollary~\ref{twist image} then implies that $\rho$ is described on generators by the formula
\[ \rho(H_{a_i}) = H_{f_\star(a_i)}^{\epsilon(f)} z^{t_i} \]
where each $t_i$ is an integer.
Since the $\sigma_i$ are all conjugate in $\an$,
we have that $t_i$ is the same for all $i$.  Conversely, any
choices of $f \in \Mod(D_{n+1})$ and $t=t_1$ determine an injective homomorphism.  Indeed, Theorem~\ref{maininj} tells us that
\[ \rho(g)Z = f (gZ)f^{-1} \]
and so the kernel of $\rho$ is contained in $Z$.  However, since $z=(\sigma_1 \cdots \sigma_n)^{n+1}$, we have
\[ \rho(z) = z \cdot z^{t(n(n+1))} = z^{1+t(n(n+1))} \]
As $n \geq 3$, we see $t(n(n+1))$ cannot be $-1$, so $\rho(z)$ is not trivial, and the kernel of $\rho$ is trivial.
Thus, we have an injection for any $t$; moreover,
the map is not surjective when $t \neq 0$: the preimage of $Z$ is contained in
$Z$, but $z \mapsto z^{1+t(n(n+1))}$, so nothing maps to $z$.

It follows that $\Aut(\an) \cong \Mod(D_{n+1})$.  This was first
proven by Dyer--Grossman \cite{dg}.  Ivanov was the first to
compute $\Aut(\an)$ from the perspective of mapping class groups
\cite{nvi}.

\p{Injections of \boldmath$\cn$.}  Again, this group has a presentation given by Figure~\ref{coxgraphpic}.  We denote the generators for $\cn$, from left to
right, by $s_1, \dots, s_n$.  The usual inclusion $\cn \to \an$ is
given by $s_1 \mapsto \sigma_1^2$ and $s_i \mapsto \sigma_i$ for
$i > 1$.

Let $\rho:\cn \to \an$ be an injective homomorphism.  There are two
punctures fixed by $\cn/Z < \Mod(S_{n+2})$.  By the moving
punctures criterion, the element $f$ given by Theorem~\ref{maininj}
must send one of these two punctures to the puncture fixed by
$\an/Z$.  Identifying $\sigma_i$ with $H_{a_i}$ as above, this means
that $f_\star$ takes $a_i$ to a 2-curve when $i >1$, and
$f_\star(a_1)$ is either a 2-curve or an $n$-curve.  By
Theorem~\ref{autcn}, $\Aut(\cn) \to \Aut(\cn/Z)$ is surjective, so we
are forced to consider these nongeometric maps $f_\star$.

As above, the homomorphism $\rho$ is given on generators by
\[ \rho(H_{a_i}) = H_{f_\star(a_i)}^{\epsilon(f)} z^{t_i} \]
Since the $s_i$ are conjugate for $i > 1$, we have that $t_i$ is the
same for these $i$; set $u=t_1$ and $v=t_2$.  Conversely, we have a
well-defined homomorphism for any $u$ and $v$.  If $g \mapsto 1$, we
again have that $g \in Z$.  Since $z=(s_1 \cdots s_n)^n$, we have $z
\mapsto z^{1+nu+n(n-1)v}$.  But there are no $u$ and $v$ which make
the latter trivial (as $n$ and $n(n-1)$ are not relatively prime), so
every choice of $u$ and $v$ leads to an injection.

\p{Injections of \boldmath$\pa$.}  We identify $\pa$ with the elements of $\an$ which fix each puncture of $D_{n+1}$.  There is a standard generating set for $\pa$, due to Artin, consisting of one Dehn twist $T_{a_{i,j}}$ for each pair of punctures of $D_{n+1}$ (see \cite{bb}).  If the punctures of $D_{n+1}$ lie in a horizontal line, then each $a_{i,j}$ can be realized as the boundary of a regular neighborhood of an arc which lies below this horizontal and connects the $i^{\mbox{\tiny th}}$ and $j^{\mbox{\tiny th}}$ punctures; note that $a_i = a_{i,i+1}$.

Let $\rho: \pa \to \an$ be an injective homomorphism.  As in the previous cases, we apply Corollary~\ref{twist image} and deduce that $\rho$ is described on generators by
\[ \rho(T_{a_{i,j}}) = T_{f_\star(a_{i,j})}^{\epsilon(f)} z^{t_{i,j}} \]
for some $f \in \Mod(S_{n+2})$.  In the case of $\pa/Z$, there are no moving punctures in $S_{n+2}$, and so the moving punctures criterion gives no restriction for the action of $f$ on the punctures of $S_{n+2}$.  We will see in Section~\ref{autsec} that in fact every $f \in \Mod(S_{n+2})$ gives rise to an automorphism of $\pa$, and so the $f$ associated to $\rho$ is arbitrary.

Conversely, since all of Artin's defining relations of $\pa$ are commutation relations (see \cite{bb}), it follows that even if the $t_{i,j}$ are all different, $\rho$ is a well-defined
homomorphism.  Again, the kernel of $\rho$ must be contained
in $Z$.  In the generators of $\pa$, $z$ can be written as
\[ (T_{a_{1,2}} T_{a_{1,3}}\cdots T_{a_{1,n+1}}) \cdots  (T_{a_{n-1,n}} T_{a_{n-1,n+1}}) (T_{a_{n,n+1}}) \]
and so we see that
\[ z \mapsto z^{1+\sum t_{i,j}} .\]
Hence, there is an affine hyperplane in $\Z^N$, where $N =
{n+1 \choose 2}$, corresponding to noninjective homomorphisms of $\pa$ into $\an$.

\p{Remark.}  The \emph{abstract commensurator} $\Comm(G)$ of a group
$G$ is the collection of isomorphisms of finite index subgroups of
$G$, where two such isomorphisms are equivalent if they agree on some
common finite index subgroup.  For a given group $\an$, $\cn$, or
$\pa$, different choices of the function $t$ give rise to distinct
elements of $\Comm(\an)$.  After the first version of this paper was
written, Leininger--Margalit proved that $\Comm(\an) \cong
\Mod(S_{n+2}) \ltimes (\Q^\times \ltimes \Q^\infty)$ \cite{lem}.


\section{Automorphisms}
\label{autsec}

In this section we construct a lift $\xi_k$ of the natural map $\Aut(L_k) \to \Aut(L_k/Z)$, where (as in the introduction) $L_k$ is the subgroup of $\an$ consisting of elements which fix $k$ particular punctures.  Using the facts that $L_1 \cong \cn$ and $L_{n+1} \cong \pa$, we will compute $\Aut(\cn)$ explicitly (reproving a result of Charney--Crisp), and we will give a generating set for $\Aut(\pa)$.

Before we begin in earnest, we note that any $L_k$ can be generated by
Artin's generators for $\pa$ plus a collection of half twists which are lifts of the elements of the symmetric group on $n+1$ letters which are in the image of $L_k$.

Charney--Crisp define a \emph{transvection} of a group $G$ with
infinite cyclic center $Z = \langle z \rangle$ to be a homomorphism $G \to G$ of the form $x \mapsto x z^{t(x)}$, where $t:G \to \Z$ is a homomorphism.  They observe that such a map is an automorphism if and only if its restriction to $\Z$ is surjective; this holds if and only if $t(z)=\pm 1$, i.e., $z \mapsto z^{\pm 1}$.  We denote by $\T(G)$ the \emph{transvection subgroup} of $\Aut(G)$.

We consider the following sequence:
\begin{equation}
\label{autseq}
1 \to \T(L_k) \to \Aut(L_k) \to \Aut(L_k/Z) \to 1
\end{equation}
To find a generating set for $\Aut(L_k)$, it suffices to construct the
lift $\xi_k:\Aut(L_k/Z) \to \Aut(L_k)$ (so the sequence is exact) and
to find generating sets for $\T(L_k)$ and $\Aut(L_k/Z)$. The group
$\T(L_k)$ can often be computed directly from a presentation of $L_k$,
and (by Theorem~\ref{main1}, say) $\Aut(L_k/Z)$ is isomorphic to the subgroup $\bar G_{k+1}$ of
$\Mod(S_{n+2})$ consisting of elements which preserve a set of $k+1$
punctures (the group $\bar G_{k+1}$ is generated by Dehn twists, half
twists, and a reflection).  In the case of $L_1 \cong \cn$, we will show that
the above exact sequence is split.

\subsection{Auxiliary groups}

Intuitively, we would like to ``blow up'' the punctures fixed by $L_k$ (the ones that are fixed by $L_k$) into boundary components so that the group $\bar G_{k+1} \cong \Aut(L_k/Z)$ cannot distinguish between the original boundary of $D_{n+1}$ and the fixed punctures.  In particular, we want $\bar G_{k+1}$ to be able to interchange $\partial D_{n+1}$ with the fixed punctures.  We now make this precise.

If one is only interested in Theorems~\ref{autcn} and~\ref{autpa}, this
subsection can be skipped; see the remarks at the end of Sections~\ref{autcnsec} and~\ref{autpasec}.

Let $\bar S$ be a sphere with $n+2$ boundary components.  We choose a set $P$ of distinguished points in $\bar S$, one in each boundary component.  Then, we define $\overline\Mod(\bar S)$ to be the group of homeomorphisms of $\bar S$ fixing $P$ as a set, modulo isotopies which fix $P$.

We fix an embedding $\bar S \to D_{n+1}$ which induces an isomorphism
on the level of curve complexes (send each boundary component to a
circle around a puncture or a circle parallel to $\partial D_{n+1}$).
We will use the same names for the curves which are equivalent under
this isomorphism (and the other isomorphisms below).

\pics{ght}{Generalized half twist.}{2}

We get the embedding $\iota: \an \to \overline\Mod(\bar S)$ as
follows.  If the generators $\sigma_i$ correspond to half twists
$H_{a_{i}}$ about the 2-curves $a_{i}$, then we define
$\iota(H_{a_{i}})$ to be the \emph{generalized half twist} about
$a_{i}$, as indicated in Figure~\ref{ghtpic}.  The generalized
half twist about a curve $a$ is denoted $\tilde H_a$.

For our definition of $\iota$ to be precise, we must specify the points of $P$.  If $\set{d_i}$ are the boundary components of $\bar S$, we choose the unique such labelling consistent with the isomorphism $\C(D_{n+1}) \cong \C(\bar S)$  and the choice of the $\set{a_{i}}$.  We draw $\bar S$ in the plane so that $d_{n+2}$ is the outer boundary component and the other $d_i$ are Euclidean circles which lie in a horizontal line.  Then, the points of $P$ are chosen to be the leftmost point of each circle (this choice is consistent with Figure~\ref{ghtpic}).

To see that $\iota$ is a homomorphism, one only needs to check the two braid relations.  The commuting relation obviously holds.  In Figure~\ref{braidpic}, we show the effect of $\iota(H_{a_{i}}H_{a_{i+1}}H_{a_{i}}) = \iota(H_{a_{i+1}}H_{a_{i}}H_{a_{i+1}})$.  We can also see that $\iota$ is injective; indeed, the map $\bar S \to D_{n+1}$ induces a left inverse $\pi: \overline\Mod(\bar S) \to \an$.  Of course, $\iota$ restricts to an injection $L_k \to \overline\Mod(\bar S)$, also called $\iota$, for any $k$.

\pics{braid}{The braid relation in $\overline\Mod(\bar S)$.}{2}

We introduce another surface $\bar S_k$, obtained by gluing punctured
disks to the $d_i$ corresponding to the punctures in
$D_{n+1}$ not fixed by $L_k$ (the surface $\bar S_k$ is a sphere with $k+1$ boundary components and $n-k+1$ punctures).  The inclusion $\bar S \to \bar S_k$ identifies
$\C(\bar S)$ with $\C(\bar S_k)$ and induces a map $\eta:\iota(L_k)
\to \overline\PMod(\bar S_k)$, where by $\overline\PMod(\bar S_k)$ we
mean the isotopy classes of homeomorphisms of $\bar S_k$ which are the
identity on the boundary.  Fixing a set of points $P \subset \partial
\bar S_k$ (one for each boundary), we can alternatively think of
$\overline\PMod(\bar S_k)$ as a normal subgroup of the group $\overline\Mod(\bar S_k)$, which consists of homeomorphisms of $\bar S_k$ fixing $P$ as a set (modulo isotopies fixing $P$).

The map $\iota_k = \eta \circ \iota$ is again injective since there is an inverse $\pi_k$ induced by gluing punctured disks to $k$ of the components of $\partial \bar S_k$.

We also want a map from $\overline\Mod(\bar S_k)$ to $\bar G_{k+1}$.  We glue punctured disks to the $k+1$ boundary components of $\bar S_k$ in order to obtain the surface $S_{n+2}$.  The inclusion of surfaces induces a surjective map $\overline\Mod(\bar S_k) \to \bar G_{k+1}$.

We encode the key relationships between all of our groups in Figure~\ref{soupgroups}.

\begin{figure}
\begin{center}
\scalebox{1}{
\xymatrix{
\an \ar@{->}[r]^\iota & \overline\Mod(\bar S) \ar@{->>}[d] \ar@{->}[r] &   \Mod(S_{n+2}) \\
L_k \ar@{->}[r] \ar@{->}[u] \ar@{->}[rd]^{\iota_k} & \overline\Mod(\bar S_k)  \ar@{->}[r] &  \bar G_{k+1} \ar@{->}[u] \\
 & \overline\PMod(\bar S_k) \ar@{->}[u] &
}
}
\caption{Groups used in the definition of $\xi_k:\Aut(L_k/Z) \to \Aut(L_k)$.}
\label{soupgroups}
\end{center}
\end{figure}

\subsection{Generalized lantern relation} In order to define our lift
$\xi_k : \Aut(L_k/Z) \to \Aut(L_k)$, we will need a relation in
$\overline\PMod(\bar S)$ called the \emph{generalized lantern relation}.
Let $\set{T_{a_{i,j}}}$ be the set of Artin generators for $\pa$, and let $\set{d_i}$ be the set of boundary components of $\bar S$.  In the language we have developed, the relation is
\begin{eqnarray*}
\iota(z) &=& \iota((T_{a_{1,2}}T_{a_{1,3}}\cdots T_{a_{1,n+1}}) \cdots
(T_{a_{n-1,n}}T_{a_{n-1,n+1}}) (T_{a_{n,n+1}})) \\
&=&  T_{d_1}^{-1}T_{d_2}^{-1} \cdots T_{d_{n+1}}^{-1} T_{d_{n+2}}
\end{eqnarray*}
(the first equality is the well-known relation in $\pa$, and the
second equality is the generalized lantern relation).  This relation appears in the work of Wajnryb \cite{bw}, who writes that this relation can be checked ``by induction (by drawing many pictures)''.  In Section~\ref{lantern}, we give a straightforward proof of the
relation.

Without reference to $\iota$, the generalized lantern relation is
simply:
\[ (T_{a_{1,2}}T_{a_{1,3}}\cdots T_{a_{1,n+1}}) \cdots
(T_{a_{n-1,n}}T_{a_{n-1,n+1}}) (T_{a_{n,n+1}}) = T_{d_1}^{n-1}T_{d_2}^{n-1} \cdots T_{d_{n+1}}^{n-1} T_{d_{n+2}}
\]
In the case of $n=2$, this relation is precisely the famous
\emph{lantern relation}, known to Dehn \cite{md}.

Stated in this alternate way, the relation exhibits an obvious asymmetry in $\pa$ between the punctures of $D_{n+1}$ and $\partial D_{n+1}$ when $n \neq 2$.  In our first description of the relation, there still is an asymmetry (in the signs), and we will see that this is what prevents us from finding a homomorphism $\Aut(\pa/Z) \to \Aut(\pa)$.

\subsection{The lift} We now define lift $\xi_k$ from $\bar G_{k+1} \cong \Aut(L_k/Z)$ to $\Aut(L_k)$.  Given an element $f \in \bar G_{k+1}$, we choose a lift $\bar f$ in $\overline\Mod(\bar S_k)$.  Since $\overline\PMod(\bar S_k)$ is normal in $\overline\Mod(\bar S_k)$, conjugation by $\bar f$ induces an automorphism $\psi_f$ of $\overline\PMod(\bar S_k)$; this automorphism is well-defined since any two lifts differ by a central element of $\overline\PMod(\bar S_k)$.  We can now define an endomorphism of $L_k$ via the composition $\pi_k \circ \psi_f \circ \iota_k$.  To see that this composition of homomorphisms is actually an automorphism of $L_k$, we will show that it is surjective.  This
suffices since $L_k$ is \emph{Hopfian}, that is, every surjective
endomorphism is an automorphism: braid groups are linear by a result of Krammer and Bigelow \cite{sb,dk} and finitely generated linear
groups are Hopfian by results of Mal'cev.

The homomorphism $\pi_k \circ \psi_f \circ \iota_k$ clearly induces a
surjection from $L_k$ to $L_k/Z$, and by the generalized lantern
relation, it also induces a surjection $Z \to Z$ ($z$ maps to either $z$ or $z^{-1}$).  It follows that
$\pi_k \circ \psi_f \circ \iota_k$ is a surjection.  We are now justified in calling the composition $\xi_k(f)$, and this defines our lifting (it is clear that $\xi_k$ is a lift).

\subsection{Generalized Artin generators} We will see below that $\Aut(L_k)$ does not preserve the conjugacy classes (in $\an$) of the generators for $L_k$.  Thus, in order to get a nice statement for how $\Aut(L_k)$ acts on $L_k$, we expand the generating set for $L_k$.

Let $a$ be either a 2-curve or an $n$-curve in $D_{n+1}$.  We denote by $g(a,\alpha)$ the following element of $L_k$:
\[ g(a,\alpha) =
\begin{cases}
H_a & a \mbox{ is a 2-curve, } H_a \in L_k \\
T_a & a \mbox{ is a 2-curve, } H_a \notin L_k \\
T_a z^{-1} & a \mbox{ is an }n\mbox{-curve with moveable puncture in exterior} \\
T_a z^{\alpha} & a \mbox{ is an }n\mbox{-curve with fixed puncture in exterior} \\
\end{cases}
\]
Each $g(a,\alpha)$ is called a \emph{generalized Artin generator} for $L_k$.  We will see that the above classification of these generators completely describes the symmetry of $L_1 \cong \cn$ and suggests an asymmetry in $L_{n+1} = \pa$.

We will need to know the image under $\iota$ of each generalized Artin generators.  We first figure this out for the usual standard generators, and then use conjugation to get the rest.

By definition $\iota(H_{a_{1}})$ is equal to $\tilde
H_{a_{1}}$. Since $T_{a_{1,2}} = H_{a_{1}}^2$, it is
straightforward to check that $\iota(T_{a_{1,2}}) = \tilde
H_{a_{1}}^2$ is equal to $T_{a_{1,2}} T_{d_1}^{-1} T_{d_2}^{-1}$
(refer to Figure~\ref{ghtpic}).  If $h \in \an$ and
$h_\star(a_{1,2})=a$, then we see that
\[
\iota(g(a,\alpha)) = \iota(h H_{a_{1,2}}^qh^{-1}) = \iota(h)\iota(H_{a_{1,2}})^q\iota(h)^{-1} = \tilde H_{h_\star(a_{1,2})}^q = \tilde H_a^q
\]
where $q \in \set{1,2}$.  The third equality holds because $\iota(h)$ and $h$ induce the same
maps of $\C(D_{n+1})$ (which is identified with $\C(\bar S$)).  As
part of our proof of the generalized lantern relation in
Section~\ref{lantern}, we will use our understanding of the action of
$\iota$ on each $H_a$ and $H_a^2$ to show that if $a$ is a curve surrounding each
puncture but the $i^{\mbox{\tiny th}}$ and $\alpha=\pm 1$, then $\iota(T_az^\alpha)$ is equal to $T_a(T_{d_i}^{-1}z)^\alpha$.

\subsection{Automorphisms of \boldmath$\cn$}\label{autcnsec}  Recall that $\cn$ is isomorphic to $L_1$, and that $\cn$ is generated by elements $s_i$ where $s_1=T_{a_{1}}$ and $s_i = H_{a_{i}}$ for $i > 1$.  We now compute the transvection subgroup of $\Aut(\cn)$ and show that our lifting $\xi_1$ from $\bar G_2 \cong \Aut(\cn/Z)$ to $\Aut(L_1) \cong \Aut(\cn)$ is a homomorphism.

In Section~\ref{cat}, we classified all transvections of $\cn$ in terms of two integers, $u$ and $v$.  These were defined by $s_1 \mapsto s_1 z^u$ and $s_2 \mapsto s_2 z^v$.  We also found that $z \mapsto z^{1+nu+n(n-1)v}$.  Again, in order for a transvection to be an automorphism, we need $z \mapsto z^{\pm 1}$.  We see that $z \mapsto z$ if and only if $nu+n(n-1)v=0$ and $z \mapsto z^{-1}$ if and only if $nu+n(n-1)v=-2$.  The latter case actually cannot happen, since $nu+n(n-1)v$ is divisible by $n \geq 3$ while $-2$ is not.  Thus, $\T(\cn) \cong \Z$.

We now want to show that $\xi_1$ is a splitting of the sequence (\ref{autseq}).  There is a homomorphism $\delta:\bar G_2 \to \Z_2$ given by the action on the two punctures of $S_{n+2}$ which are fixed setwise by $G_2$.  Recall that $\epsilon:\bar G_2 \to \Z_2$ is the homomorphism which records whether or not elements are orientation preserving.

Given $f \in \bar G_2 \cong \Aut(L_1/Z)$, we have the following simple formula for the action of $\xi_1(f)$ on the generalized Artin generators of $L_1$:
\[ \xi_1(f)(g(a,\alpha)) = g(f_\star(a),\delta(f)\alpha)^{\epsilon(f)}. \]
It follows easily that $\xi_1$ is a homomorphism.  Checking the above formula is straightforward for each of the 4 types of generators; the work was done in their careful classification into the 4 types.  To give the idea, we check the formula for the fourth type of generator (this is the only case where $\delta$ is important).  Let $a$ be a curve which surrounds all punctures but the fixed one.
\begin{eqnarray*}
\iota_1(g(a,\alpha)) = \iota_1(T_az^{\alpha}) &=& T_{a} (T_{d_1}^{-1} T_{d_{n+2}})^\alpha \\
\psi_f \circ \iota_1(g_a) &=& (T_{f_\star(a)} (T_{f_\star(d_1)}^{-1} T_{f_\star(d_{n+2})})^\alpha)^{\epsilon(f)} \\
\xi_1(f)(T_az^\alpha) = \pi_1 \circ \psi_f \circ \iota_1(T_az^\alpha) &=&
\begin{cases} (T_{f_\star(a)} z^{\alpha})^{\epsilon(f)} & \mbox{ if } \delta(f) = +1
\\  (T_{f_\star(a)} z^{-\alpha})^{\epsilon(f)}  &\mbox{ if } \delta(f) = -1\end{cases} \\
&=& g(f_\star(a),\delta(f)\alpha)^{\epsilon(f)}
\end{eqnarray*}

We now have that $\Aut(\cn) \cong \bar G_2 \ltimes \Z$.  The group $\bar G_2$ is isomorphic to $(\Z_2 \times \Z_2) \ltimes G_2$, where, as in the introduction, the group $G_2$ is the group of orientation preserving elements of $\bar G_2$ which fix two particular punctures.  Thus, we can write $\Aut(\cn)$ as $((\Z_2 \times \Z_2) \ltimes G_2) \ltimes Z$.  As noted by Charney--Crisp, the elements of $G_2$ commute with the transvections of $\Aut(\cn)$, and so, finally, we obtain Theorem~\ref{autcn}: 
$\Aut(\cn) \cong (\Z_2 \times \Z_2) \ltimes (G_2 \times \Z)$.

\p{Remark.} In the special case of $k=1$, we can give a more
straightforward definition of the lift $\xi_k$.  Given an $f \in \bar
G_2$, we define $\xi_1(f)$ directly by the formula
$\xi_1(f)(g(a,\alpha)) = g(f_\star(a),\delta(f)\alpha)^{\epsilon(f)}$.
Using the presentation of $\cn$, and the generalized lantern relation,
one can directly check that this defines a homomorphism $\bar G_2 \to
\Aut(\cn)$.

\subsection{Automorphisms of \boldmath$\pa$}
\label{autpasec}
 Having constructed $\xi_{n+1}$ we have completed our proof of Theorem~\ref{autpa}.  We now address the question of whether or not there is a splitting $\Aut(\pa/Z) \to \Aut(\pa)$.  More specifically, we will explain why the map $\xi_{n+1}$ is not a splitting.  Since $\xi_{n+1}$ seems like the most natural candidate for a splitting, we conjecture that there is no splitting.

Let $a$ be a curve in $D_{n+1}$ which surrounds all punctures except the first.  Let $g \in \Mod(S_{n+2}) \cong \Aut(\pa/Z)$ be an element whose lift $\bar g \in \overline\Mod(\bar S)$ satisfies $\bar g(d_1)=d_1$ and $\bar g(d_{n+2}) = d_2$.  Similarly to our calculations for $\cn$, we can check that $\xi_{n+1}(g)$ takes $T_az^{-1}$ to $T_{g_\star(a)}$.  Let $f \in \Mod(S_{n+2})$ be such that $f \circ g(a) = a$, but $\bar f \circ \bar g(d_1) = d_{n+2}$ and $\bar f \circ \bar g(d_{n+2}) = d_1$.  Then $\xi_{n+1}(f)$ takes $T_{g_\star(a)}$ to $T_az^{-1}$.  However, $\xi_{n+1}(fg)$ takes $T_az^{-1}$ to $T_az$.  Thus, $\xi_{n+1}$ is not a homomorphism.

\p{Remark.} There is a simpler proof that $\Aut(\pa)$ surjects
onto $\Aut(\pa/Z)$.  Given any $f \in \Mod(S_{n+2}) \cong
\Aut(\pa/Z)$, one can use the presentation of $\pa$ to directly show
that for any choice of $\set{t_{i,j}}$, the map given by $T_{a_{i,j}}
\mapsto T_{f_\star(a_{i,j})} z^{t_{i,j}}$ is a homomorphism (a
convenient presentation to use for this is the ``modified'' Artin
presentation in \cite{mm}).  Thus, in order to obtain a lift of $f$,
one only needs to choose $\set{t_{i,j}}$ so that $z \mapsto z^{\pm
  1}$; that is, $\sum t_{i,j}$ needs to be 0 or $-2$.  Since we don't
have explicit presentations for the other $L_k$, this method does not work in general.

We also note here that an earlier version of this paper contained an incorrect computation of $\Aut(\pa)$.


\section{Generalized lantern relation}
\label{lantern}

We now prove the ``generalized lantern relation'', used in the proof
of Theorem~\ref{autpa}.  We freely use the notation of
Section~\ref{autsec}.

Our goal is to understand $\iota(z)$.  We think of $z$ as a product of
elements
\[ g_i = T_{a_{i,i+1}}T_{a_{i,i+2}}\cdots T_{a_{i,n+1}} \]
and we will first understand each $\iota(g_i)$ individually.  We draw
$\bar S$ in the plane as in Section~\ref{autsec}.  This allows us
to see the $a_{i,j}$ in $\bar S$ exactly as they appear in $D_{n+1}$.

\pics{zigzag}{}{3}

We can think of each $T_{a_{i,j}}$ as a ``push map'', where the
$i^{\mbox{\tiny th}}$ boundary component moves around the
$i^{\mbox{\tiny th}}$ boundary component, while travelling clockwise
inside $a_{i,j}$ in such a way that it never turns (Figure~\ref{ghtpic} represents the halfway point of this push map).

\pics{easier}{}{3}

We can thus think of $g_i$ as a product of these push maps (see
Figure~\ref{zigzagpic}).  We now see the following intuitive relation:
$g_i$ can also be obtained by pushing the $i^{\mbox{\tiny th}}$
boundary component around the $n-i+1$ boundary components to its right
all at once (see Figure~\ref{easierpic}).  We then observe that this
latter push map is equivalent to
\[ T_{c_{i,n+1}} T_{c_{i+1,n+1}}^{-1} T_{d_i}^{-1} \]
where the curves are as shown in Figure~\ref{twistspic}.  This intuitive relation
(which is already an interesting relation in the mapping class group)
is explained more formally in the remark below (see also \cite{mm}).

\pics{twists}{}{3}

We can now compute $\iota(z)$ as the product of the $\iota(g_i)$:
\[ (T_{c_{1,n}} T_{c_{2,n}}^{-1} T_{d_1}^{-1}) (T_{c_{2,n}}
T_{c_{3,n}}^{-1}T_{d_2}^{-1}) \cdots (T_{c_{n,n+1}} T_{c_{n,n}}^{-1}
T_{d_{n+1}}^{-1}) \] All of the $T_{c_{i,j}}$ elements cancel except
the first, which is equal to $z$, and the last, which is equal to
$T_{d_{n}}^{-1}$.  Thus, $\iota(z)$ is equal to the product of $z$
with $T_{d_1}^{-1} \cdots T_{d_{n+1}}^{-1}$, and this is exactly the
generalized lantern relation.

We notice that, applying the map $\pi$ to $\iota(z)$, we see that we
have proven that the product of the $g_i$ is indeed equal to $z$ in $\pa$.

Also, since $\iota$ takes conjugates of the $T_{a_{i,j}}$ to the corresponding conjugates of the $\iota(T_{a_{i,j}})$, the same holds for the conjugates of the $g_i$.  This fact was used in Section~\ref{autsec}.

\p{Remark.}  We now explain a more formal framework for proving the
intuitive relation that pushing a disk around two loops is the same as
pushing it around a composite loop.  Let $X$ be the
subset of the configuration space of $n+1$ ordered points in the unit
tangent bundle of the disk, where each point lies in a different
fiber.  There is a natural map from $\pi_1(X) \to \overline\Mod(\bar S)$ (the
projection to the disk of each point in $X$ specifies the location of a particular
boundary component and the vector specifies the rotation; if we like,
we can replace boundary components with rigid disks in the sphere).  The
relation described above simply follows from the fact that this map is
a homomorphism.  Putman has observed that relations in $\pi_1(X)$ thus
give rise to many different ``generalized lantern relations''.


\section{Superinjective maps of curve complexes}

Let $S_m$ be a sphere with $m \geq 5$ punctures, and let $\phi$ be a
superinjective map of $\C(S_m)$.  We will prove Theorem~\ref{superinj}, i.e., that $\phi$ is induced by a unique element of $\Mod(S_m)$.  The basic strategy is to show that $\phi$ preserves the topological types of curves (Lemma~\ref{3c}), that $\phi$ preserves the simplest type of nontrivial intersection between curves (Lemma~\ref{adj}), and then to use induction to show that $\phi$ is surjective (Proposition~\ref{surjectivity}).  Theorem~\ref{superinj} is then a consequence of the following theorem of Korkmaz \cite{mk}:

\begin{thm}
Let $m \geq 5$.  Every automorphism of $\C(S_m)$ is induced by a unique element of $\Mod(S_m)$.
\end{thm}

A \emph{side} of a curve $a$ in $S_m$ is one of the two connected components of $S_m-a$.  The curve $a$ is called a \emph{k-curve} if the minimum of the numbers of punctures on each side is $k$. Two 2-curves $a$ and $b$ in $S_m$ are said to be {\em adjacent} if $i(a,b)=2$.

\begin{lem}[Sides]\label{sides}
If $a$ and $b$ are two distinct curves which lie on the
same side of a curve $w$, then $\phi(a)$ and $\phi(b)$ lie on the same
side of $\phi(w)$.\end{lem}

\bpf

Choose a curve $d$ which intersects $a$ and $b$, but not $w$.  Since
$\phi$ is superinjective, $\phi(d)$ intersects $\phi(a)$ and
$\phi(b)$, but not $\phi(w)$, and so the lemma follows.
\epf

\begin{lem}[2-curves]\label{2c}
If $a$ is a 2-curve, then $\phi(a)$ is a 2-curve.
\end{lem}

\bpf

Choose a pants decomposition $\set{a=c_1, c_2, \dots, c_{m-3}}$.
Applying Lemma~\ref{sides}, we see that $\phi(c_2), \dots, \phi(c_{m-3})$
must all lie on the same side of $\phi(a)$.  It follows that $\phi(a)$ is
a 2-curve.
\epf

\begin{lem}[k-curves]\label{3c}
If $w$ is a k-curve, then $\phi(w)$ is a k-curve.
\end{lem}

\bpf

By Lemma~\ref{2c}, we may assume that $w$ is not a 2-curve.  Any pants decomposition $\P$ containing $w$ has $k-2$ curves on one side and $m-k-2$ curves on the other side.  By Lemma~\ref{sides} and the injectivity of $\phi$, the curve $\phi(w)$ must either be a $k$-curve or a 2-curve.  Thus, it suffices to rule out the latter possibility.

First, suppose that $w$ has an even number of punctures on one of its sides. Choose maximal collections $C$ and $C'$ of disjoint 2-curves on each side of $w$.  The set $C \cup C'$, and hence $\phi(C) \cup \phi(C')$, realizes the maximal number of disjoint 2-curves in $S_m$.  Since $\phi$ is injective, $\phi(w)$ cannot be a 2-curve.

Now suppose that $w$ has an odd number of punctures on both of its
sides. There are two cases.  In the first case, one of the sides of
$w$ has at least five punctures.  Choose maximal collections $C_0$ and
$C_1$ of 2-curves on each side of $w$, where $C_0$ has at least two
2-curves. If $\phi(w)$ is a 2-curve, then the collection $C = C_0 \cup
C_1 \cup \{w\}$ maps to a maximal collection of disjoint
2-curves. Choose $c \in C_0$ and an adjacent 2-curve $x$ with
$i(x,y)=0$ for all $y \in C - \set{c}$. Since $\phi(c)$ and $\phi(x)$ are distinct 2-curves and $\phi(C)$ is maximal, $\phi(c) \cup \phi(x)$ must separate two other curves of $\phi(C)$. But for any curves $c_1,c_2 \in C$ not equal to $c$, we can find a curve $y$ which intersects $c_1$ and $c_2$ but does not intersect $c$ or $x$ (such a curve will not exist in the case that $w$ has 3 punctures on each side). Mapping these curves forward, we contradict the fact that $\phi(x) \cup \phi(c)$ separates.

In the second case, both sides of $w$ have 3 punctures and $S_m = S_6$. Choose 2-curves $c_0$, $x_0$, $c_1$, and $x_1$, each disjoint from $w$, so that $c_0$ and $x_0$ are adjacent and lie on one side of $w$ and so that $c_1$ and $x_1$ are adjacent and lie on the other side of $w$.  If $\phi(w)$ is a 2-curve, then $\phi(c_i) \cup \phi(x_i)$ separates $\phi(w)$ and $\phi(c_{i-1})$ (indices read modulo 2).  We see that $\phi(c_0)$ and $\phi(c_1)$ intersect, contradicting the fact that $c_0$ and $c_1$ are disjoint.
\epf

\begin{lem}[Sides II]
\label{sides 2}
Let $w$ be a curve in $S_m$, let $F_1$ and $F_2$ be the sides of $w$,
and let $F_1'$ and $F_2'$ be the sides of $\phi(w)$.  Up to
renumbering, $F_i$ is homeomorphic to $F_i'$ and $\phi$ induces maps
from the vertices of $\C(F_i)$ to those of $\C(F_i')$ for $i=1,2$.
\end{lem}

\bpf

By Lemma~\ref{3c}, we have that $F_i$ is homeomorphic to $F_i'$ for $i=1,2$.  If $\P$ is a pants decomposition containing $w$, then $\P-\{w\}$ restricts to pants decompositions of $F_1$ and $F_2$, and an application of Lemma~\ref{sides} completes the proof.
\epf

\begin{lem}[Adjacency]\label{adj}
If $a$ and $b$ are adjacent, then $\phi(a)$ and $\phi(b)$ also adjacent.\end{lem}

\bpf We claim that 2-curves $a$ and $b$ are adjacent
if and only if there exists a curve $w$ and curves $x$ and $y$ so
that: $a$ and $b$ lie on a thrice-punctured side of $w$, $x$
intersects $a$ and $w$ but not $b$ and not $y$, and $y$ intersects $b$
and $w$ but not $a$ and not $x$.  By Lemma~\ref{3c}, Lemma~\ref{sides 2} and the definition
of superinjectivity, all of these properties are preserved by $\phi$,
and thus the lemma will follow.

One direction is easy: if $a$ and $b$ are adjacent, then we can find
curves $w$, $x$, and $y$ which satisfy the given properties.  Now
suppose that there exist curves $w$, $x$, and $y$ which satisfy the
given properties.  We restrict our attention to the side of $w$
containing $a$ and $b$.  On this subsurface $S'$, $x$ and $y$ are
collections of disjoint arcs.  Note that on a thrice-punctured disk,
there can be at most three families of disjoint parallel arcs.
However, since $a$ is a curve disjoint from $y$, arcs of $y$ can only
appear in one of these families.  The same is true for $x$, and we see
that the arcs of $x$ are not parallel to the arcs of $y$.  Thus, $a$
must lie in the component of $S'-y$ which is a twice punctured disk.
There is only one such curve.  Likewise, there is only one choice for
$b$, and we see that $a$ and $b$ are adjacent.
\epf

We alter the definition of the curve complex for $S_4$ so it is the
flag complex with vertices for curves in $S_4$ and edges corresponding
to adjacency.  It is well-known that $\C(S_4)$ is isomorphic to the
classical \emph{Farey graph} (see \cite{ynm}).  Further it is not hard to see that an injective simplicial map of the Farey graph to itself is determined by what it does to a single triangle.  We thus have the following fact.

\begin{lem}
\label{complexity 1}
Any injective simplicial map $\C(S_4) \to \C(S_4)$ is surjective.
\end{lem}

As discussed at the start of this section, the following proposition completes the proof of Theorem~\ref{superinj}.

\begin{prop}
\label{surjectivity}
$\phi$ is surjective.
\end{prop}

\bpf

We proceed by induction on $m$, starting with base case $m=4$, which is covered by Lemma~\ref{complexity 1}.  Now assume that $m \geq 5$, and that the proposition is true for all spheres with fewer punctures.

Let $c$ be any curve in $S_m$.  By Lemma~\ref{sides 2}, $\phi$ induces superinjective maps from the curve complexes of the sides of $c$ to the curve complexes of the sides of $S_m-\phi(c)$; if any of the components are homeomorphic to $S_4$, then we apply Lemma~\ref{complexity 1}, and if any of the components are homeomorphic to $S_3$, then we simply throw it out.

By induction, we conclude that these induced maps are surjective.  In
other words, for any curve in the image of $\phi$, the entire set of
points joined to that curve by an edge in $\C(S_m)$ is in the image of
$\phi$.  Since $\C(S_m)$ is connected, it follows that $\phi$ is
surjective.
\epf

\p{Remark.} An earlier version of this paper used a different argument between
Lemma~\ref{adj} and the conclusion of the proof.  The approach,
already used in several papers \cite{ni,mk,ei,ei2,ei3,bem}, was to
show that $\phi$ induces a map of the arc complex and to find an element of $\Mod(S_m)$ which agrees with $\phi$ on some maximal simplex of the
arc complex (injective simplicial maps of the arc complex are
determined by their action on a single maximal simplex).  We did this by quoting parts of
Korkmaz's proof that automorphisms of $\C(S_m)$ induce injective
simplicial maps of the arc complex.

The idea of showing the surjectivity of $\phi$ and directly applying
the theorem of Korkmaz came much later, and gives a much simpler way
of completing the proof.  Concurrently with the revision of this
paper, the induction argument was used by the authors to give a very
short argument that all superinjective maps are surjective \cite{si}.

We also remark that the proof can also be simplified using the idea of
the \emph{adjacency graph}, introduced by Behrstock--Margalit and
Shackleton \cite{bem} \cite{kjs}.  However, so as not to overly confuse the
chronological relationships between these papers, we refrain from
employing this useful tool here.


\bibliographystyle{plain}

\bibliography{coartin}

\end{document}